%% file: sundial.tex
\documentclass[a4paper]{llncs}
\usepackage{url, color, graphicx, amsmath, colordvi, subfigure, epsfig}
\usepackage{capt-of}
\usepackage{verbatim}
\usepackage{algorithm}
\usepackage[noend]{algpseudocode}

\begin{document}
\title{Sundial: Using Sunlight to Reconstruct Global Timestamps}
\author{Jayant Gupchup$^\ast$  R\u{a}zvan Mus\u aloiu-E.$^\ast$ Alex Szalay$^\ddagger$ Andreas Terzis$^\ast$}
\institute{
Computer Science Department$^\ast$ \quad Physics and Astronomy Department$^\ddagger$ \\
Johns Hopkins University \\ 
{\texttt{\{gupchup,razvanm,terzis\}@jhu.edu}$^\ast$ \quad
  \texttt{szalay@jhu.edu}$^\ddagger$ \\
}
}

\maketitle

\begin{abstract}
This paper investigates postmortem timestamp reconstruction in
environmental monitoring networks. In the absence of a
time-synchronization protocol, these networks use multiple pairs of
(local, global) timestamps to retroactively estimate the motes' clock
drift and offset and thus reconstruct the measurement time series. We
present {\em Sundial}, a novel offline algorithm for reconstructing
global timestamps that is robust to unreliable global clock
sources. Sundial reconstructs timestamps by correlating annual solar
patterns with measurements provided by the motes' inexpensive light
sensors. The surprising ability to accurately estimate the length of
day using light intensity measurements enables Sundial to be robust to
arbitrary mote clock restarts. Experimental results, based on multiple
environmental network deployments spanning a period of over 2.5 years,
show that Sundial achieves accuracy as high as 10 parts per million
(ppm), using solar radiation readings recorded at 20 minute intervals.
\end{abstract}

\input{intro}

\input{problem}

\input{prot}

\input{eval}

\input{relwork}

\input{concl}

\input{ack}

\begin{small}

\end{small}

\end{document}

%% file: intro.tex
\section{Introduction}
\label{sec:intro}

A number of environmental monitoring applications have demonstrated
the ability to capture environmental data at scientifically-relevant
spatial and temporal scales~\cite{LUYFJournal,TPS+05}. These
applications do not need online clock synchronization and in the
interest of simplicity and efficiency often do not employ one. Indeed,
motes do not keep any global time information, but instead, use their
local clocks to generate local timestamps for their measurements.
Then, a postmortem timestamp reconstruction algorithm retroactively
uses (local, global) timestamp pairs, recorded for each mote
throughout the deployment, to reconstruct global timestamps for all
the recorded local timestamps. This scheme relies on the assumptions
that a mote's local clock increases monotonically and the global clock
source (e.g., the base-station's clock) is completely
reliable. However, we have encountered multiple cases in which these
assumptions are violated. Motes often reboot due to electrical shorts
caused by harsh environments and their clocks restart. Furthermore,
basestations' clocks can be desynchronized due to human and other
errors. Finally the basestation might fail while the network continues
to collect data.

We present {\em Sundial}, a robust offline time reconstruction
mechanism that operates in the absence of any global clock source and
tolerates random mote clock restarts. Sundial's main contribution is a
novel approach to reconstruct the global timestamps using only the
repeated occurrences of day, night and noon.  We expect Sundial to
work alongside existing postmortem timestamp reconstruction
algorithms, in situations where the basestations' clock becomes
inaccurate, motes disconnect from the network, or the basestation
fails entirely. While these situations are infrequent, we have
observed them in practice and therefore warrant a solution. We
evaluate Sundial using data from two long-term environmental
monitoring deployments. Our results show that Sundial reconstructs
timestamps with an accuracy of one minute for deployments that are
well over a year.

%% file: problem.tex
\section{Problem Description}
\label{sec:prob}

\begin{figure}[t]
  \centering
  \includegraphics[scale=0.5]{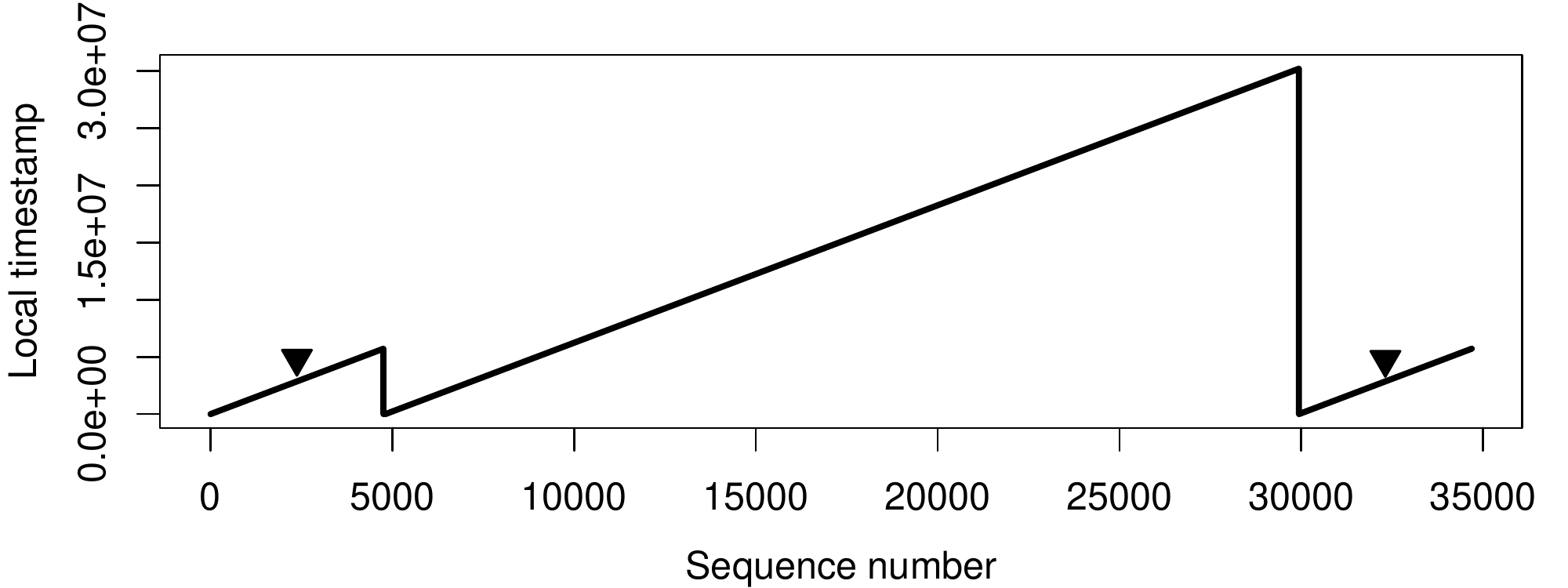}
  \caption {An illustration of mote reboots, indicated by clock
    resets. Arrows indicate the segments for which anchor points are
    collected.} 
  \label{fig:reboot}
\end{figure}  

The problem of reconstructing global timestamps from local timestamps
applies to a wide range of sensor network applications that correlate
data from different motes and external data sources. This problem is
related to mote clock synchronization, in which motes' clocks are
persistently synchronized to a global clock source. However,
In this work, we focus
on environmental monitoring applications that do not use online time
synchronization, but rather employ postmortem timestamp reconstruction
to recover global timestamps.

\subsection{Recovering Global Timestamps}
\label{subsec:mtr}

As mentioned before, each mote records measurements using its local
clock which is not synchronized to a global time source. During the
lifetime of a mote, a basestation equipped with a global clock
collects multiple pairs of (local, global) timestamps. We refer to
these pairs as \emph{anchor points}\footnote{We ignore the
  transmission and propagation delays associated with the anchor point
  sampling process.}. Furthermore, we refer to the series of local
timestamps as $LTS$ and the series of global timestamps as $GTS$. The
basestation maintains a list of anchor points for each mote and is
responsible for reconstructing the global timestamps using the anchor
points and the local timestamps.

The mapping between local clock and global clock can be described by
the linear relation $GTS = \alpha \cdot LTS + \beta$, where $\alpha$
represents the slope and $\beta$ represents the intercept (start
time). The basestation computes the correct $\alpha$ and $\beta$ for
each mote using the anchor points. Note that these $\alpha$ and
$\beta$ values hold, if and only if the mote does not reboot. In the
subsections that follow, we describe the challenges encountered in
real deployments where the estimation of $\alpha$ and $\beta$ becomes
non-trivial.

\begin{figure}[t]
  \centering
  \includegraphics[scale=0.5]{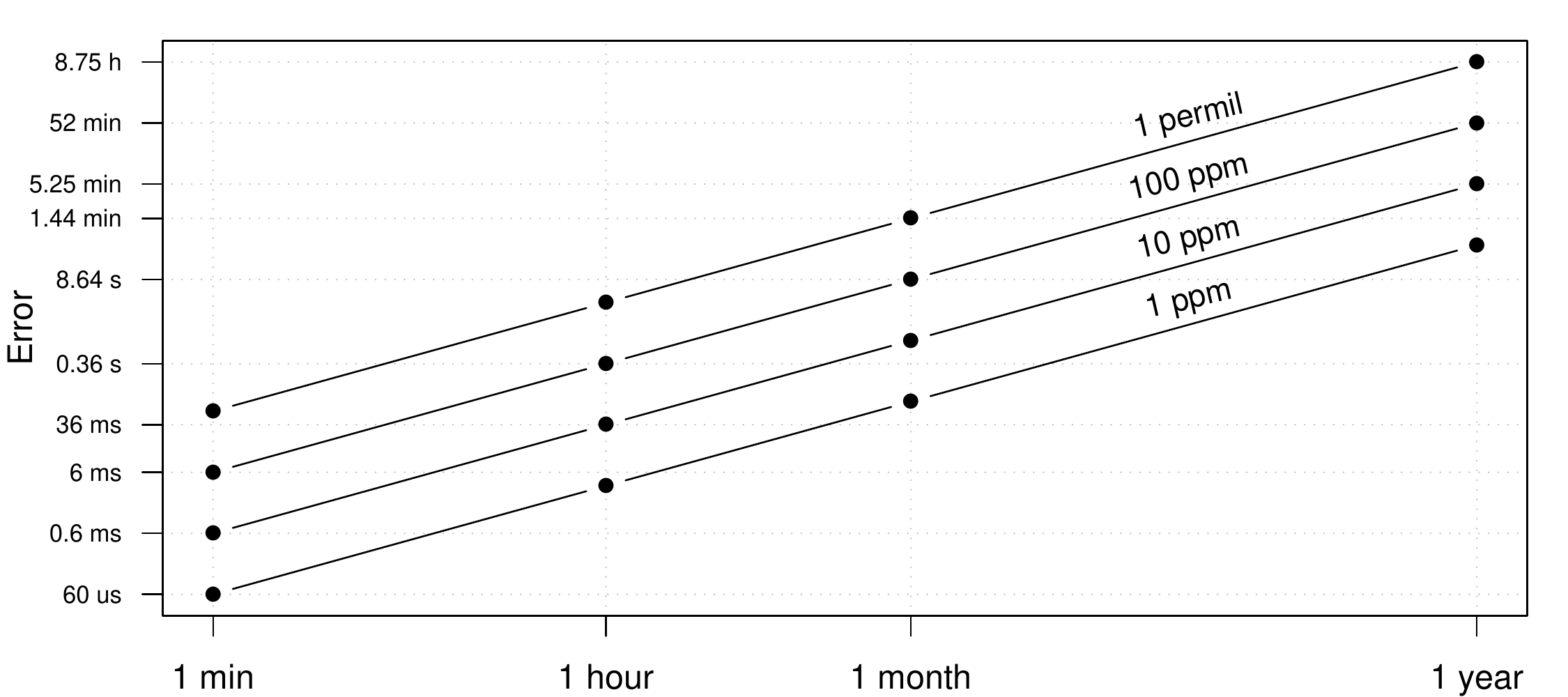}
  \caption{Time reconstruction error due to $\alpha$ estimation errors 
    as a function of the deployment lifetime.}
  \label{fig:drift}
\end{figure}

\subsection{Problems in Timestamp Reconstruction}
\label{subsec:reboot}

The methodology sketched in Section \ref{subsec:mtr} reconstructs the
timestamps for blocks of measurements where the local clock increase
monotonically. We refer to such blocks as \emph{segments}. Under ideal
conditions, a single segment includes all the mote's measurements.
However, software faults and electrical shorts (caused by moisture in
the mote enclosures) are two common causes for unattended mote reboots.
The mote's local clock resets after a reboot and when this happens we
say that the mote has started a new segment.

When a new segment starts, $\alpha$ and $\beta$ must be
recomputed. This implies that the reconstruction mechanism described
above must obtain at least two anchor points for each
segment. However, as node reboots can happen at arbitrary times,
collecting two anchor points per segment is not always
possible. Figure~\ref{fig:reboot} shows an example where no anchor
points are taken for the biggest segment, making the reconstruction of
timestamps for that segment problematic. In some cases we found that
nodes rebooted repeatedly and did not come back up immediately. Having
a reboot counter helps recover the segment chronology but does not
provide the precise start time of the new segment.

Furthermore, the basestation is responsible for providing the global
timestamps used in the anchor points. Our experience shows that
assuming the veracity of the basestation clock can be
precarious. Inaccurate basestation clocks can corrupt anchor points
and lead to bad estimates of $\alpha$ and $\beta$ introducing errors
in timestamp reconstruction. Long deployment exacerbate these
problems, as Figure \ref{fig:drift} illustrates: an $\alpha$ error of
100 parts per million (ppm) can lead to a reconstruction error of 52
minutes over the course of a year.

\subsection{A Test Case}
\label{subsec:testcase}

\begin{figure}[t]
  \centering
  \includegraphics[scale=0.5]{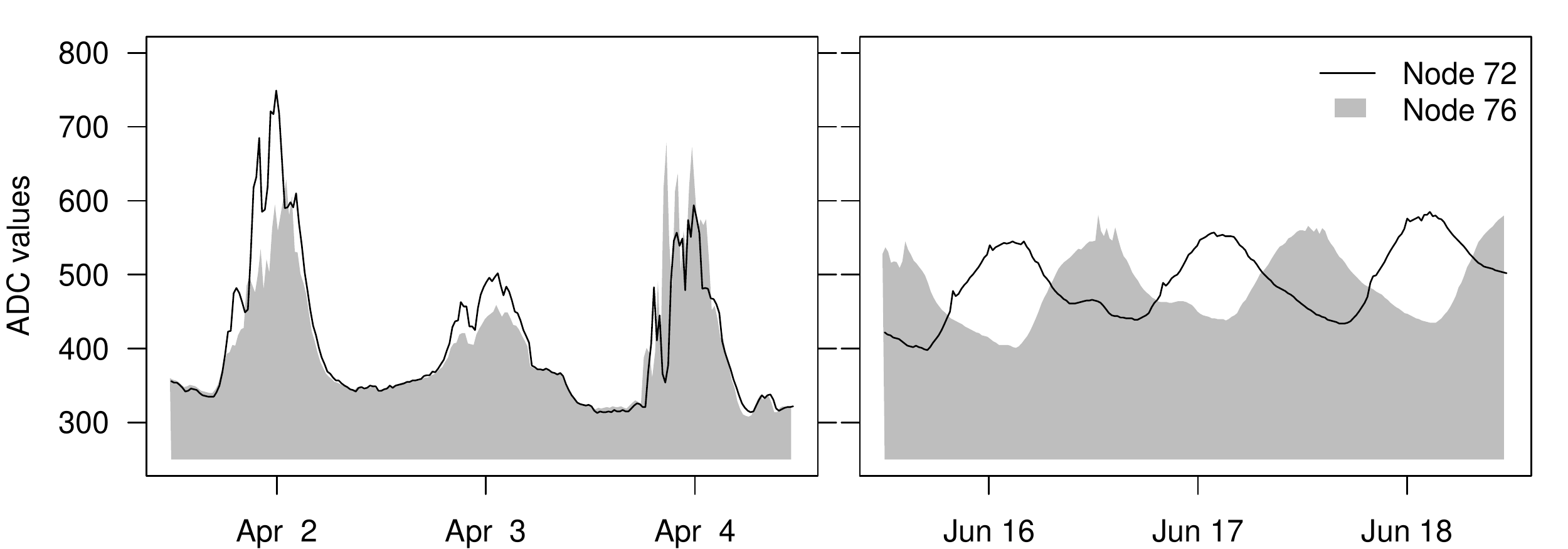}
  \caption {Ambient temperature data from two motes from the {\em L}
    deployment. The correlation of temperature readings in the left
    panel indicates consistent timestamps at the segment's
    start. After two months, the mote's reading become inconsistent 
    due to inaccurate $\alpha$ estimates.}
  \label{fig:julia}
\end{figure}

Our {\em Leakin Park} deployment (referred to as {\em ``L''}) provides
an interesting case study of the problems described above. The {\em L}
deployment comprised six motes deployed in an urban forest to study the
spatial and temporal heterogeneity in a typical urban soil
ecosystem. The deployment spanned over a year and a half, providing us
with half a million measurements from five sensing modalities.
We downloaded data from the sensor nodes very infrequently using a
laptop PC and collected anchor points only during these downloads. One
of the soil scientists in our group discovered that the ambient
temperature values did not correlate among the different
motes. Furthermore, correlating the ambient temperature with an
independent weather station, we found that the reconstruction of
timestamps had a major error in it.  

Figure \ref{fig:julia} shows data from two ambient temperature sensors
that were part of the {\em L} deployment.  Node 72 and 76 show
coherence for the period in April, but data from June are completely
out-of-sync.  We traced the problem back to the laptop acting as the
global clock source. We made the mistake of not synchronizing its
clock using NTP before going to the field to download the data. As a
result the laptop's clock was off by 10 hours, giving rise to large
errors in our $\alpha$ and $\beta$ estimates and thereby introducing
large errors in the reconstructed timestamps.  To complicate matters
further, we discovered that some of the motes had rebooted a few times
between two consecutive downloads and we did not have any anchor
points for those segments of data.

%% file: prot.tex

\begin{algorithm}[t]
\caption{Robust Global Timestamp Reconstruction (RGTR)}\label{clockrecon-algo}
\begin{algorithmic}\scriptsize

\Statex \textbf{constants}

\State $Q$ \Comment{Constant used to identify anchor points for the segment}
\State $\delta_{HIGH},\ \delta_{LOW},\ \delta_{DEC}$ \Comment{Constants used in iterative fit}

\Statex

\Procedure{ClockFit}{$ap$}
  \State $(r, i) \gets (0, 0)$
  \State $q \gets$ \Call{HoughQuantize}{$ap$}
  \For{\textbf{each} $\gamma$ \textbf{in} \Call{Keys}{$q$}}
    \State $s \gets$ \Call{Size}{$q\{\gamma\}$}
    \If {$s > r$}
      \State $(r, i) \gets (s, \gamma)$
    \EndIf
  \EndFor
  \State \textbf{return} \Call{ComputeAlphaBeta}{$q\{i\}$}
\EndProcedure

\Statex

\Procedure{HoughQuantize}{$ap$}
  \State $q \gets \{\}$  \Comment {Map of empty sets}
  \For{\textbf{each} $(lts_i, gts_i)$ \textbf{in} $ap$}
    \For{\textbf{each} $(lts_j, gts_j)$ \textbf{in} $ap$ \textbf{and} $(lts_j, gts_j) \ne (lts_i, gts_i)$}
    	\State $\alpha \gets (gts_j-gts_i)/(lts_j-lts_i)$
    	\If {$0.9 \le \alpha \le 1.1$}   \Comment {Check if part of the same segment}
    		\State $\beta \gets gts_j - \alpha \cdot lts_j$
    		\State $\gamma \gets ROUND (\beta/Q)$  
		\State 	\Call{Insert}{$q\{\gamma\}, (lts_i, gts_i)$}
		\State 	\Call{Insert}{$q\{\gamma\}, (lts_j, gts_j)$}
    	\EndIf		
    \EndFor
  \EndFor
  \State \textbf{return} $q$
\EndProcedure

\Statex

\Procedure{ComputeAlphaBeta}{$ap$}
   \State $\delta \gets \delta_{HIGH}$
   \State $bad \gets \{\}$
   \While {$\delta > \delta_{LOW}$}
	\State $(\alpha, \beta) \gets $ \Call{LLSE}{$ap$}
	\For{\textbf{each} $(lts, gts) \in ap$ \textbf{and} $(lts, gts) \notin bad$}
		\State $residual \gets (\alpha \cdot lts+\beta)-gts$
		\If {$residual \ge |\delta|$}
			\State \Call{Insert}{$bad, (lts, gts)$}
		\EndIf
	\EndFor
	\State $\delta \gets \delta - \delta_{DEC}$
   \EndWhile
   \State \textbf {return} $(\alpha, \beta$)
\EndProcedure
\end{algorithmic}
\label{rgrtalgo}
\end{algorithm}


\section{Solution}
\label{sec:prot}

The test case above served as the motivation for a novel methodology
that robustly reconstructs global timestamps. The Robust Global
Timestamp Reconstruction (RGTR) algorithm, presented in
Section~\ref{subsec:rgtr}, outlines a procedure to obtain robust
estimates of $\alpha$ and $\beta$ using anchor points that are
potentially unreliable. We address situations in which the basestation
fails to collect any anchor points for a segment through a novel
method that uses solar information alone to generate anchor points. We
refer to this mechanism as Sundial.

\subsection{Robust Global Timestamp Reconstruction (RGTR)}
\label{subsec:rgtr}

Having a large number of anchor points ensures immunity from inaccurate
ones, provided they are detected. Algorithm~\ref{rgrtalgo}
describes the Robust Global Timestamp Reconstruction (RGTR) algorithm
that achieves this goal. RGTR takes as input a set of anchor points
($ap$) for a given segment and identifies the anchor points that
belong to that segment, while censoring the bad ones. Finally, the
algorithm returns the $(\alpha, \beta)$ values for the segment. RGTR
assumes the availability of two procedures: \textsc{Insert} and
\textsc{Llse}.  The \textsc{Insert$(x,y)$} procedure adds a new
element, $y$, to the set $x$. The Linear Least Square
Estimation~\cite{Duda}, \textsc{Llse} procedure takes as input a set
of anchor points belonging to the same segment and outputs the
parameters $(\alpha, \beta)$ that minimize the sum of square
errors.

RGTR begins by identifying the anchor points for the segment.  The
procedure \textsc{HoughQuantize} implements a well known feature
extraction method, known as the Hough Transform~\cite{Duda-Hough}. The
central idea of this method is that anchor points that belong to the
same segment should fall on a straight line having a slope of $\sim
1.0$. Also, if we consider pairs of anchors (two at a time) and
quantize the intercepts, anchors belonging to the same segment should
all collapse to the same quantized value (bin). \textsc{HoughQuantize}
returns a map, $q$, which stores the anchor points that collapse to
the same quantized value. The key (stored in $i$) that contains the
maximum number of elements contains the anchor points for the segment.

\begin{figure}[t]
\begin{minipage}[t]{0.45\textwidth}
  \centering
  \includegraphics[width=\textwidth]{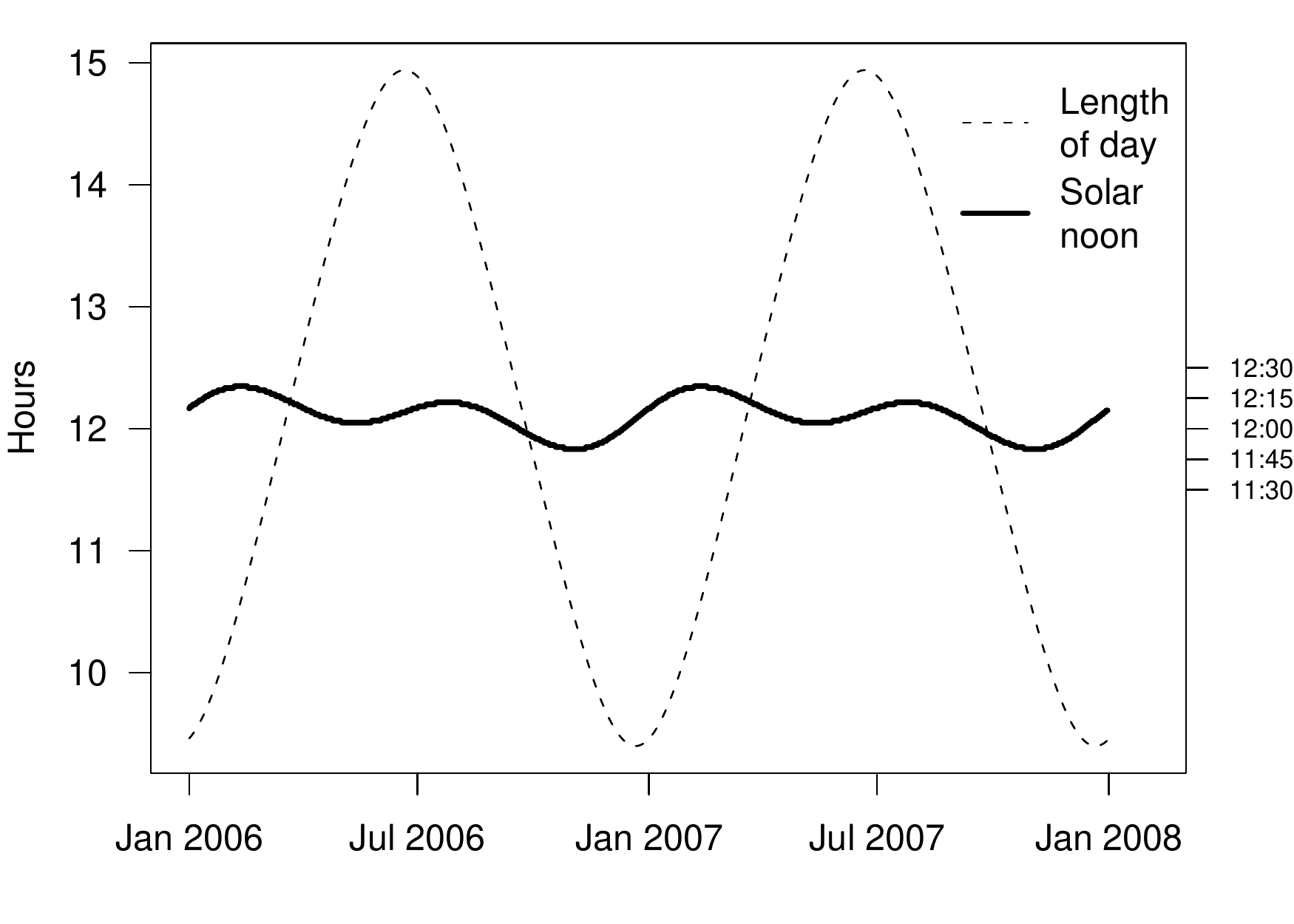}
  \caption{The solar (model) length of day (LOD) and noon pattern 
for a period of two years for the latitude of our deployments.}
  \label{sfig:model}
\end{minipage}
\hfill
\begin{minipage}[t]{0.45\textwidth}
  \centering
  \includegraphics[width=\textwidth]{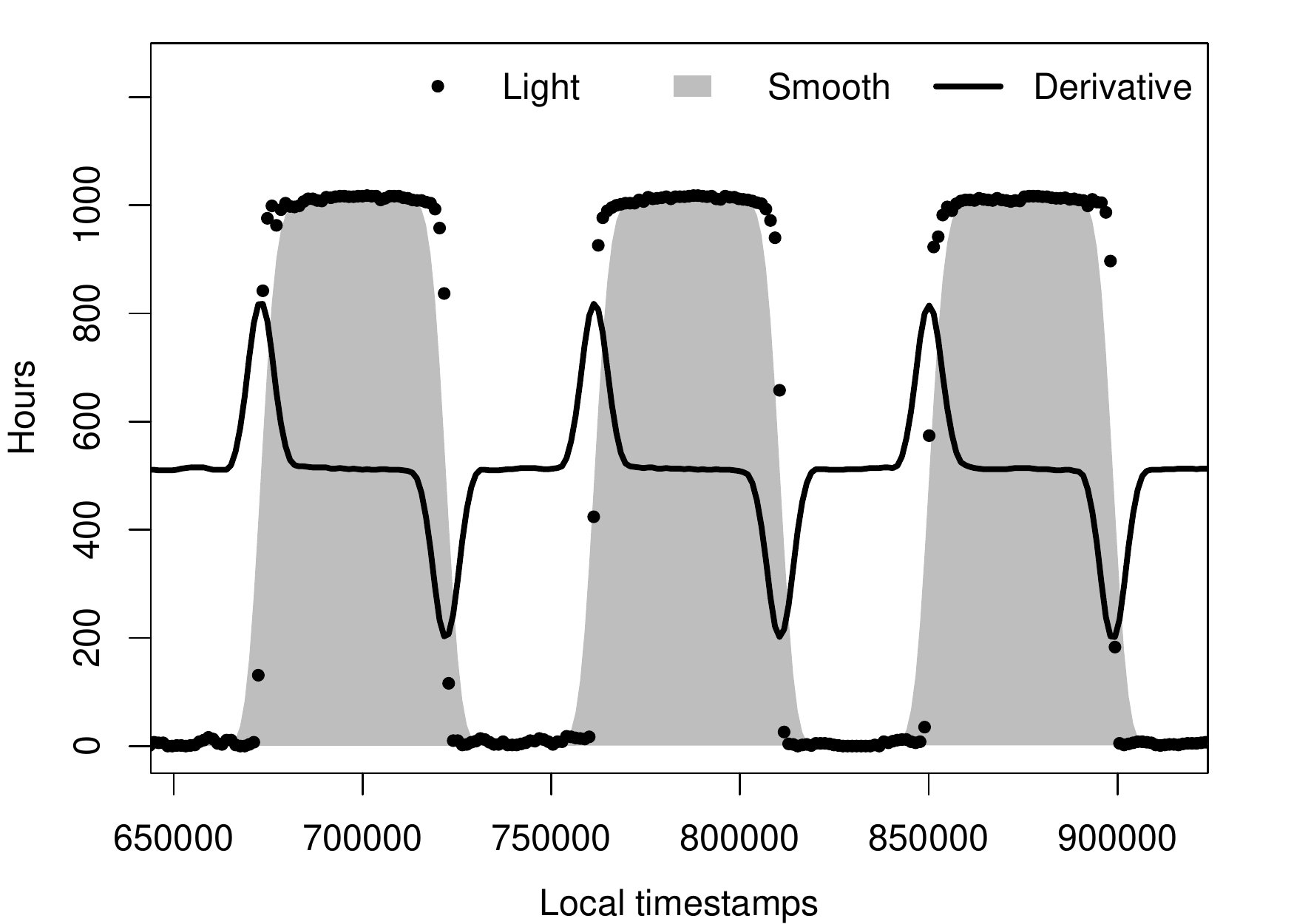}
  \caption{The light time series (raw and smoothed) and its
    first derivative. The inflection points represent sunrise and
    sunset.} 
  \label{sfig:der}
\end{minipage}
\end{figure}

Next, we invoke the procedure \textsc{ComputeAlphaBeta} to compute
robust estimates of $\alpha$ and $\beta$ for a given segment. We begin
by creating an empty set, $bad$.  The set $bad$ maintains a list of
all anchor points that are detected as being outliers and do not
participate in the parameter estimation. This procedure is iterative
and begins by estimating the fit ($\alpha, \beta$) using all the
anchor points.  Next, we look at the residual of all anchor points
with the fit. Anchor points whose residuals exceed the current
threshold, $\delta$, are added to the $bad$ set and are excluded in
the next iteration fit.  Initially, $\delta$ is set conservatively to
$\delta_{HIGH}$.  At the end of every iteration, the $\delta$
threshold is lowered and the process repeats until no new entries are
added to the $bad$ set, or $\delta$ reaches $\delta_{LOW}$.

\subsection{Sundial}
\label{subsec:sundial}

The parameters of the solar cycle (sunrise, sunset, noon) follow a
well defined pattern for locations on Earth with a given
latitude. This pattern is evident in Figure~\ref{sfig:model} that
presents the length of day (LOD) and solar noon for the period between
January 2006 and June 2008 for the latitude of the {\em L}
deployment. Note that the LOD signal is periodic and
sinusoidal. Furthermore, the frequency of the solar noon signal is
twice the frequency of the LOD signal. We refer the reader
to~\cite{SOLAR} for more details on how the length of day can be
computed for a given location and day of the year.

The paragraphs that follow explain how information extracted from
our light sensors can be correlated with known solar information to
reconstruct the measurement timestamps.

\begin{figure}[t]
\begin{center}
  \includegraphics[scale=0.5]{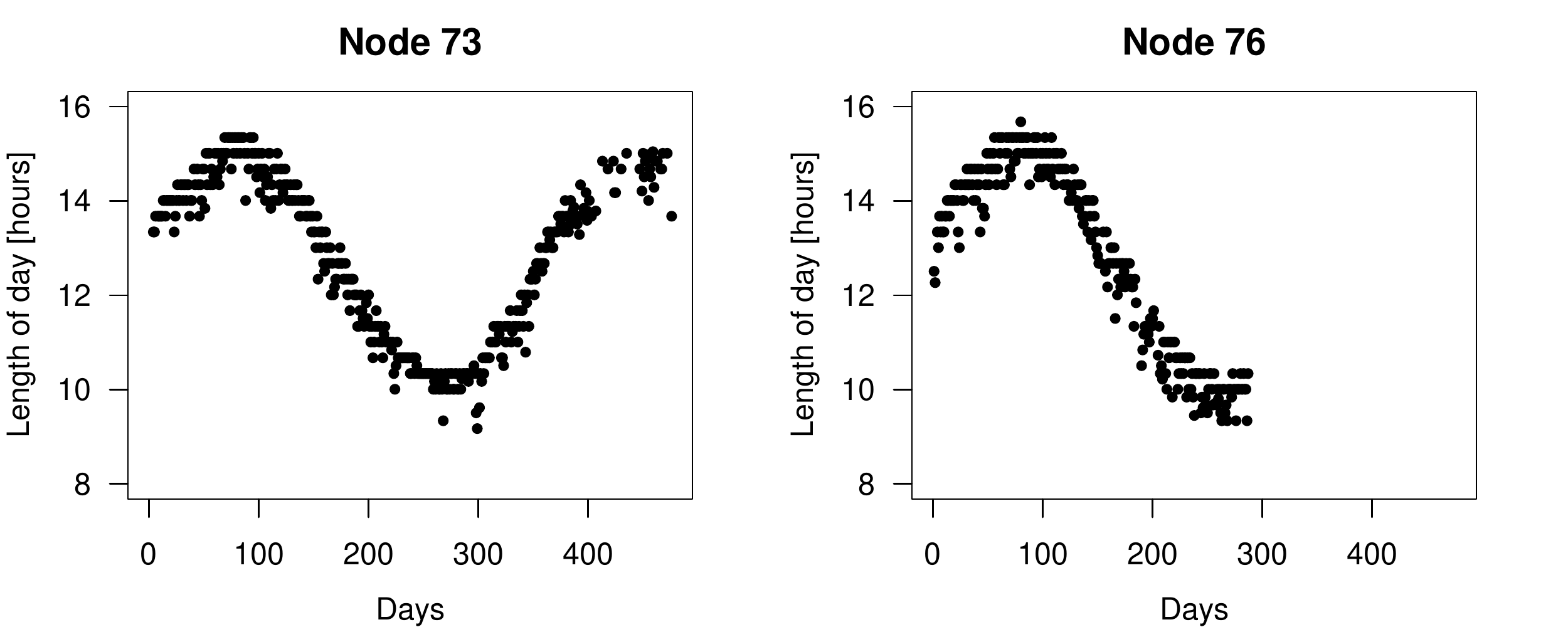}
  \caption {The length of day pattern for two long segments belonging
    to different nodes. Day 0 represents the start-time for each of
    the segments.  }
  \label{fig:lod}
\end{center}
\end{figure}

\subsubsection{Extracting light patterns: }

We begin by looking at the time series $L_i$ of light sensor readings
for node $i$. $L_i$ is defined for a single segment in terms of the
local clock. First, we create a smooth version of this series, to
remove noise and sharp transients. Then, we compute the first
derivative for the smoothed $L_i$ series, generating the $D_i$
time-series. Figure~\ref{sfig:der} provides an illustration of a
typical $D_i$ series overlaid on the light sensor series ($L_i$). One
can notice the pattern of inflection points representing sunrise and
sunset.  The regions where the derivative is high represent mornings,
while the regions where the derivative is low represent evenings. For
this method, we select sunrise to be the point at which the derivative
is maximum and sunset the point at which the derivative is
minimum. Then, LOD is given as the difference between sunrise and
sunset, while noon is set to the midpoint between sunrise and sunset.

The method described above accurately detects noon time. However, the
method introduces a constant offset in LOD detection and it
underestimates LOD due to a late sunrise detection and an early sunset
detection. The noon time is unaffected due to these equal but opposite
biases. In practice, we found that a simple thresholding scheme works
best for finding the sunrise and sunset times. The light sensors'
sensitivity to changes simplifies the process of selecting the
appropriate threshold. In the end, we used a hybrid approach whereby
we obtain noon times from the method that uses derivatives and LOD
times from the thresholding method. The net result of this procedure
is a set of noon times and LOD for each day from the
segment's start in terms of the local clock.  Figure~\ref{fig:lod}
shows the LOD values obtained for two different node segments after
extracting the light patterns.

\begin{figure}[t]
  \centering
  \includegraphics[scale=0.5]{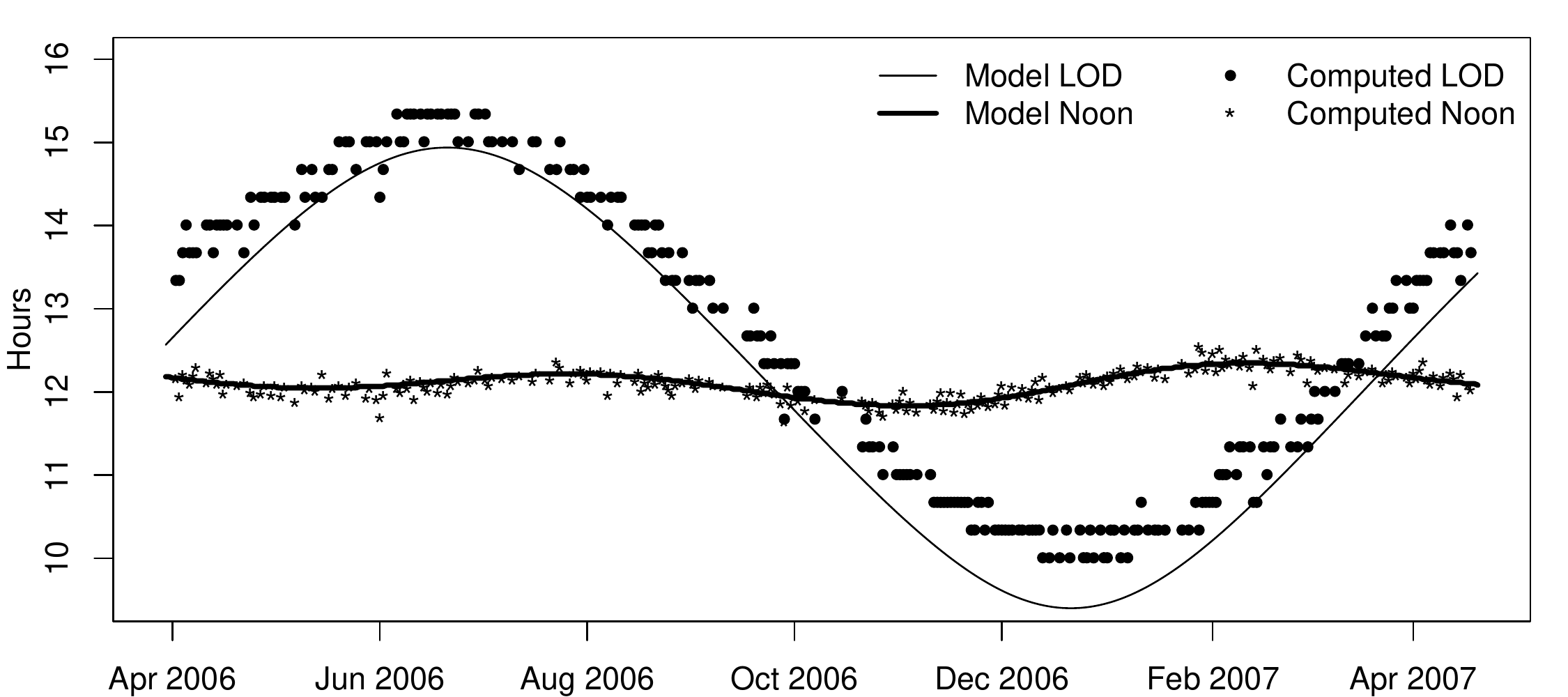}
  \caption{An illustration of the computed LOD and noon values for the lag with
    maximum correlation with the solar model.}
  \label{fig:validation}
\end{figure}

\subsubsection{Solar reconstruction of clocks: }

The solar model provides the LOD and noon values in terms of the
global clock ($LOD_{GT}$), while the procedure described in the
previous paragraph extracts the LOD and noon values from light sensor
measurements in terms of the motes' local clocks ($LOD_{LT}$).  In
order to find the best possible day alignment, we look at the
correlation between the two LOD signals ($LOD_{GT}, LOD_{LT}$) as a
function of the lag (shift in days). The lag that gives us the maximum
correlation ($\rho_{max}$) is an estimate of the day alignment.
Mathematically, the day alignment estimate (lag) is obtained as

$$ \arg\max_{lag} \text{Cor}(LOD_{GT}, LOD_{LT}, lag) $$

\noindent where $\text{Cor}(X,Y,s)$ is the correlation between time
series $X$ and $Y$ shifted by $s$ time
units. Figure~\ref{fig:validation} presents an example of the match
between model and computed LOD and noon times achieved by the lag with
the highest correlation. The computed LOD time series tracks the one
given by the solar model. One also observes a constant shift between
the two LOD patterns, which can be attributed to the horizon
effect. For some days, canopy cover and weather patterns cause the
extracted LOD to be underestimated.  However, as the day alignment is
obtained by performing a cross-correlation with the model LOD pattern,
the result is robust to constant shifts. Furthermore,
Figure~\ref{fig:validation} shows that the equal and opposite effect
of sunrise and sunset detection ensures that the noon estimation in
unaffected in the average case.

After obtaining the day alignment, we use the noon information to
generate anchor points.  Specifically, for each day of the segment we
have available to us the noon time in local clock (from the light
sensors) and noon time in global clock (using the model). RGTR can
then be used to obtain robust values of $\alpha$ and $\beta$. This fit
is used to reconstruct the global timestamps. As
Figure~\ref{sfig:model} suggests, the noon times change slowly over
consecutive days as they oscillate around 12:00.  Thus, even if the
day estimate is inaccurate, due to the small difference in noon times,
the $\alpha$ estimate remains largely unaffected. This implies that
even if the day alignment is not optimal, the time reconstruction
within the day will be accurate, provided that the noon times are
accurately aligned.  The result of an inaccurate lag estimate is that
$\beta$ is off by a value equal to the difference between the actual
day and our estimate. In other words, $\beta$ is off by an integral
and constant number of days (without any skew) over the course of the
whole deployment period.

\begin{figure}[t]
  \centering
  \includegraphics[height=3.3in]{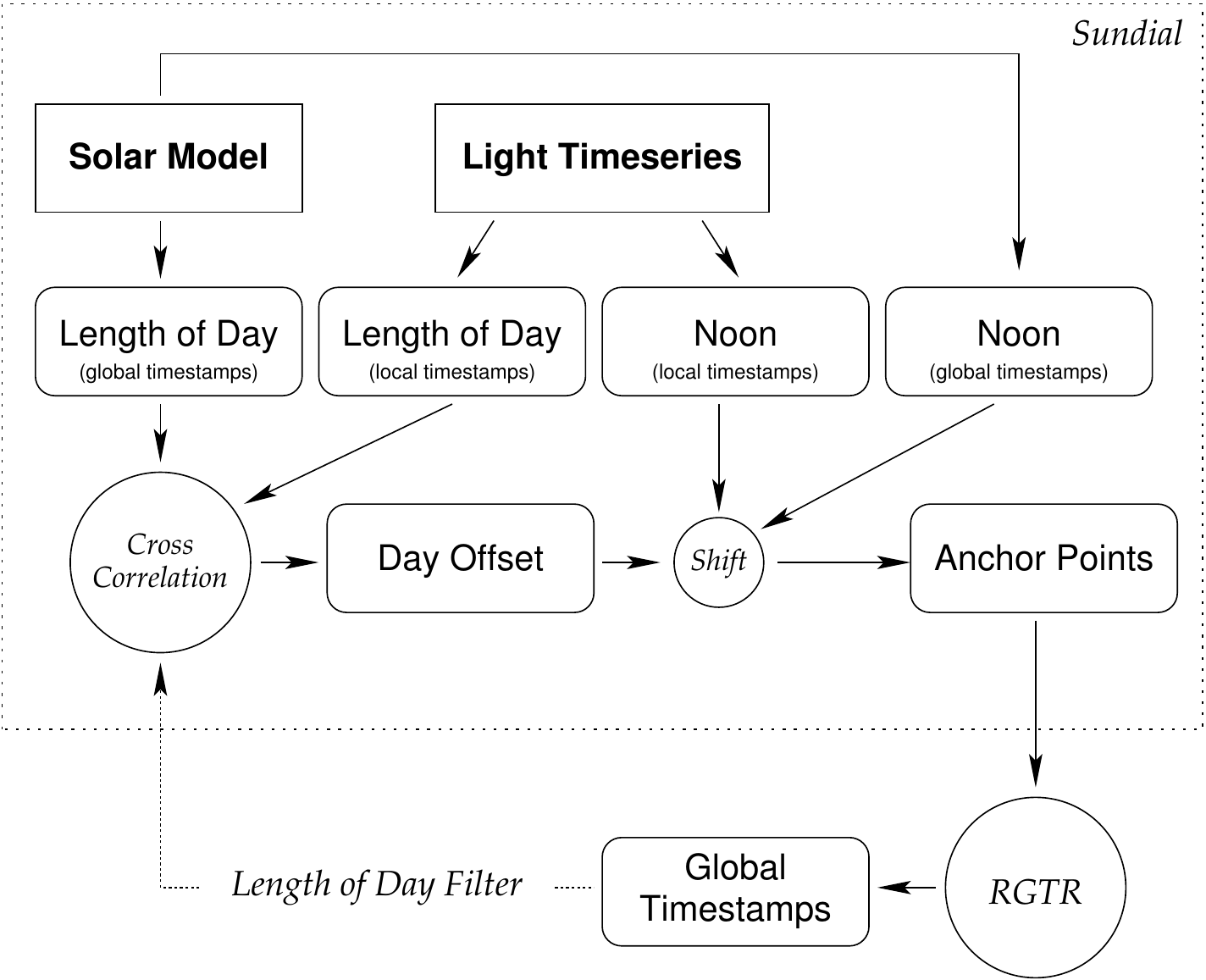}
  \caption{The steps involved in reconstructing global timestamps
    using Sundial.}
  \label{fig:flow}
\end{figure}

We find that this methodology is well suited in finding the correct
$\alpha$. To improve the $\beta$ estimate, we perform an iterative
procedure which works as follows. For each iteration, we obtain the
best estimate fit ($\alpha, \beta$). We convert the motes' local
timestamps into global timestamps using this fit. We then look at the
difference between the actual LOD (given by the model) and the current
estimate for that day. If the difference between the expected LOD and
the estimate LOD exceeds a threshold, we label that day as an outlier.
We remove these outliers and perform the LOD cross-correlation to
obtain the day shift (lag) again. If the new lag differs from the lag
in the previous iteration, a new fit is obtained by shifting the noon
times by an amount proportional to the new lag.  We iterate until the
lag does not change from the previous iteration. Figure \ref{fig:flow}
shows a schematic of the steps involved in reconstructing global
timestamps for a segment.

\begin{figure}[t]
\begin{center}
  \includegraphics[scale=0.5]{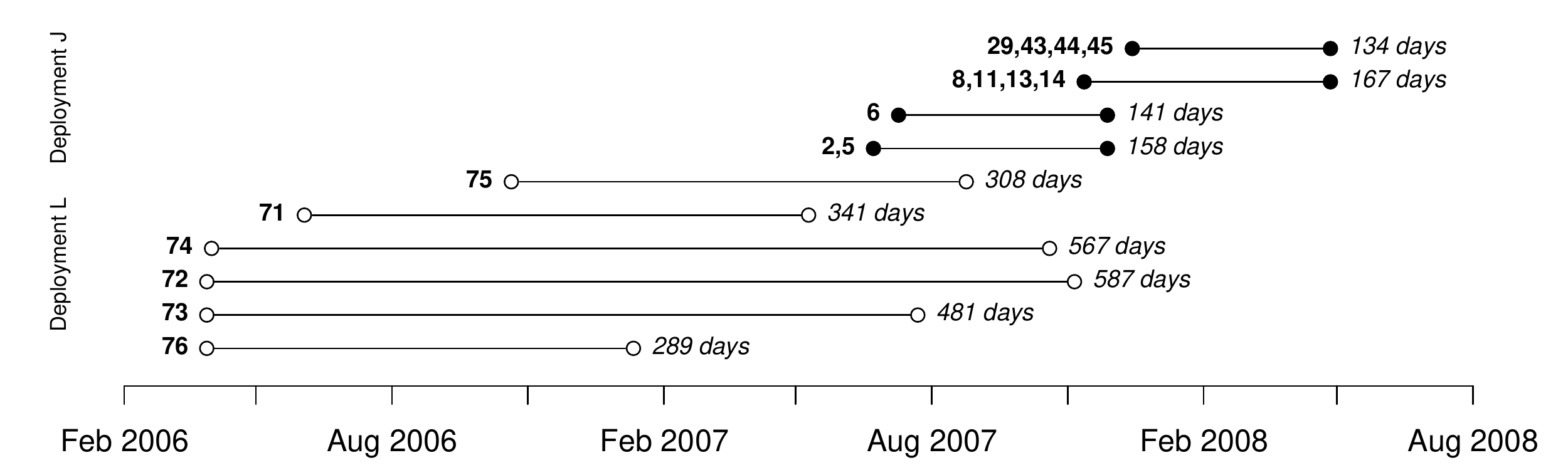}
  \caption{Node identifiers, segments and length of each
    segment (in days) for the two deployments used in the evaluation.} 
  \label{fig:deploy}
\end{center}
\end{figure}

%% file: eval.tex
\section{Evaluation}
\label{sec:eval}

We evaluate the proposed methodology using data from two deployments.
Deployment {\em J} was done at the Jug Bay wetlands sanctuary along
the Patuxent river in Anne Arundel County, Maryland. The data it
collected is used to study the nesting conditions of the Eastern Box
turtle ({\em Terrapene carolina}) \cite{AGU07}. Each of the motes was
deployed next to a turtle nest, whereas some of them have a clear view
of the sky while others are under multiple layers of tree
canopy. Deployment {\em L}, from Leakin Park, is described in
Section~\ref{subsec:testcase}.

Figure~\ref{fig:deploy} summarizes the node identifiers, segments, and
segment lengths in days for each of the two deployments. Recall that a
segment is defined as a block of data for which the mote's clock
increases monotonically. Data obtained from the {\em L} dataset
contained some segments lasting well over 500 days. The {\em L}
deployment uses MicaZ motes~\cite{micazspecs}, while the {\em J}
deployment uses TelosB motes~\cite{PSC05}.  Motes 2, 5, and 6
from Deployment {\em J} collected samples every 10 minutes. All
other motes for both deployments had a sampling interval of 20
minutes. In addition to its on-board light, temperature, and humidity
sensors, each mote was connected to two soil moisture and two soil
temperature sensors.

In order to evaluate Sundial's accuracy, we must compare the
reconstructed global timestamps it produces, with timestamps that
are known to be accurate and precise. Thus we begin our evaluation by
establishing the ground truth.

\begin{figure}[t]
\begin{minipage}[t]{0.45\textwidth}
  \centering
  \includegraphics[width=\textwidth]{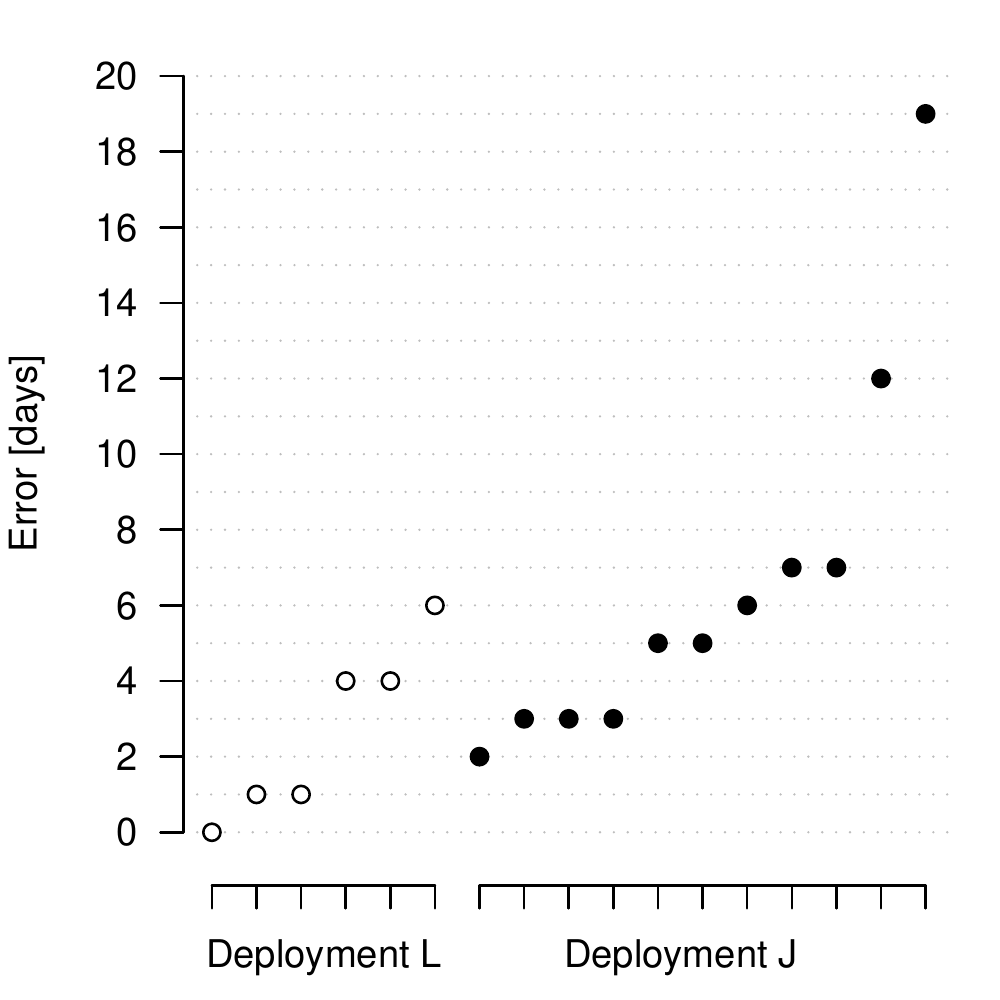}
  \caption{Error in days for different motes from the {\em L} and {\em
      J} deployments.}
  \label{sfig:dayerr}
\end{minipage}
\hfill
\begin{minipage}[t]{0.45\textwidth}
  \centering
  \includegraphics[width=\textwidth]{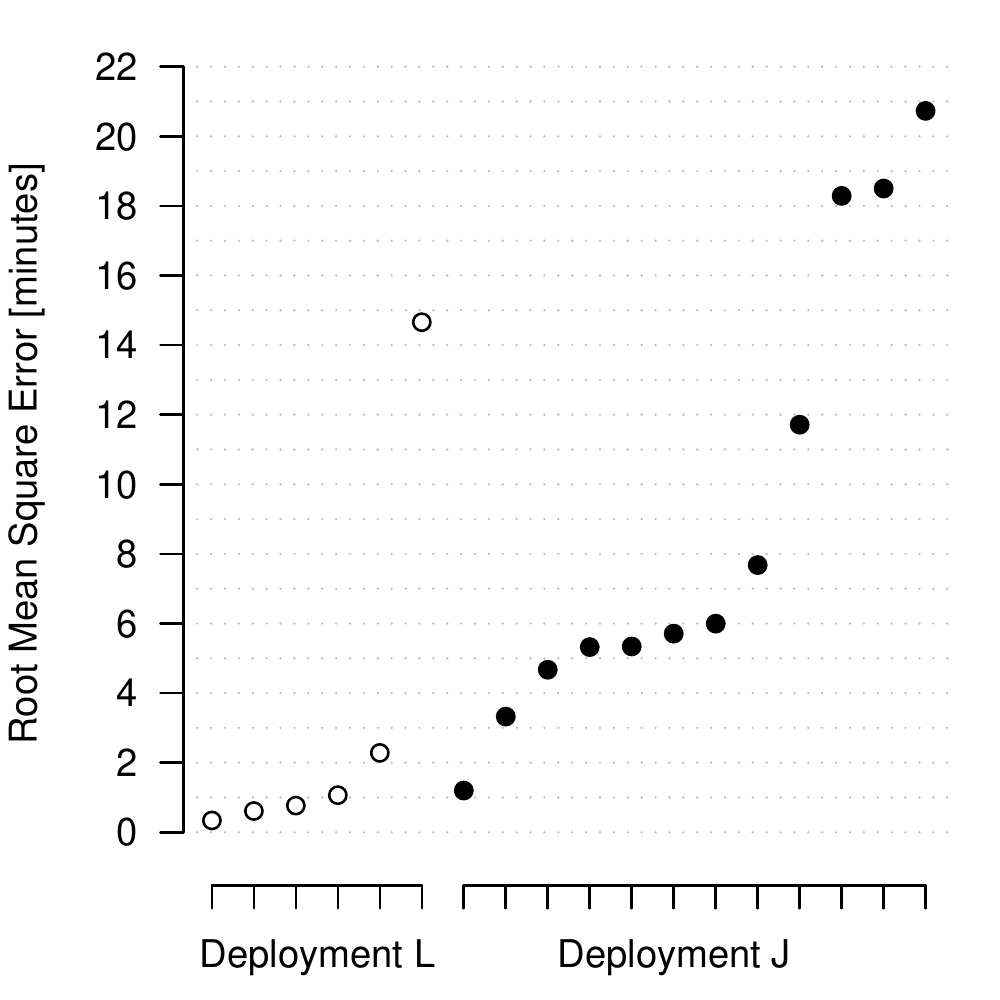}
  \caption{Root mean square error in minutes $(RMSE_{min})$.}
  \label{sfig:minerr}
\end{minipage}
\end{figure}

\subsection{Ground Truth}
	
For each of the segments shown in Figure~\ref{fig:deploy}, a set of
good anchor points (sampled using the basestation) were used to obtain
a fit that maps the local timestamps to the global timestamps. We
refer to this fit as the {\em Ground truth fit}. This fit was
validated in two ways. First, we correlated the ambient temperature
readings among different sensors. We also correlated the motes'
measurements with the air temperature measurements recorded by nearby
weather stations. The weather station for the {\em L} deployment was
located approximately 17 km away from the deployment site~\cite{BWI},
while the one for the {\em J} deployment was located less than one km
away~\cite{JBPR}. Considering the proximity of the two weather
stations we expect that their readings are strongly correlated to the
motes' measurements.

Note that even if the absolute temperature measurements differ, the
diurnal temperature patterns should exhibit the same behavior thus
leading to high correlation values. Visual inspection of the
temperature data confirmed this intuition.  Finally, we note that due
to the large length of the segments we consider, any inconsistencies
in the ground truth fit would become apparent for reasons similar to
the ones provided in Section~\ref{subsec:reboot}.

\begin{figure}[t]
\begin{minipage}[t]{0.45\textwidth}
  \centering
  \includegraphics[width=\textwidth]{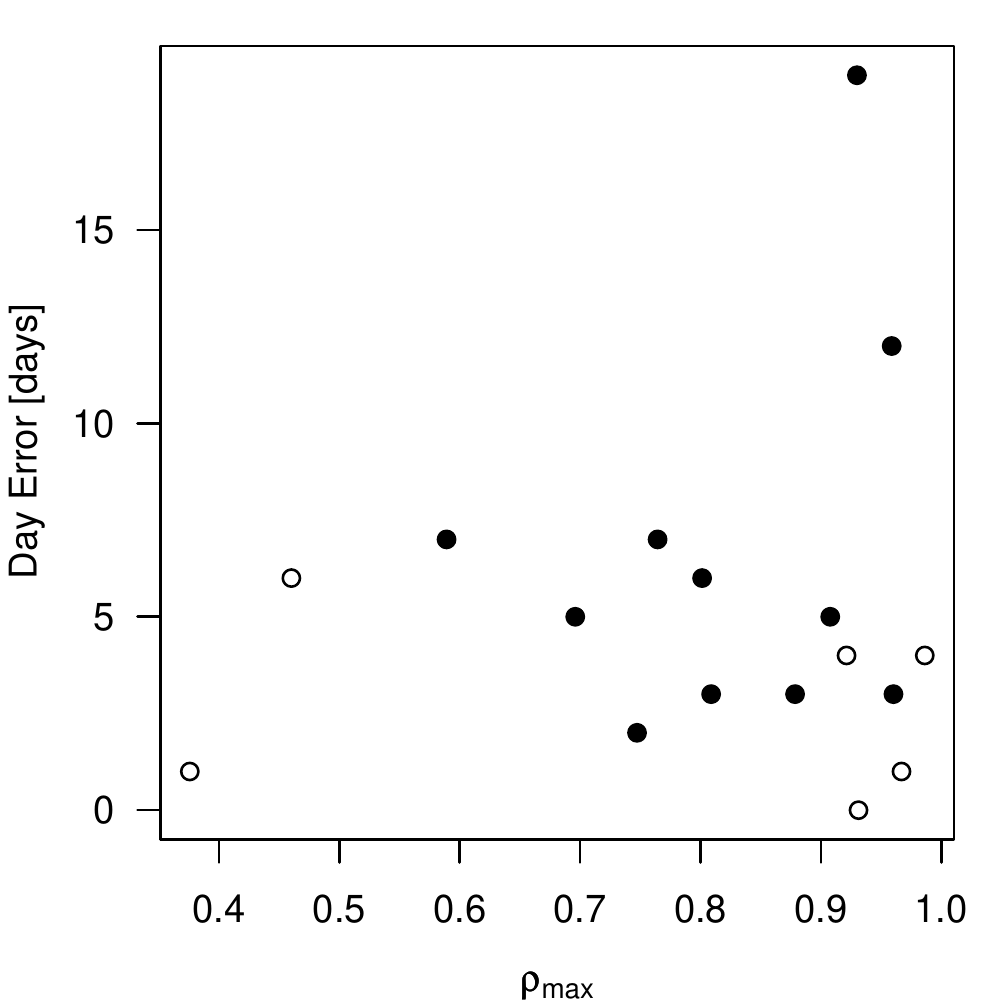}
  \caption{Relation between $\rho_{max}$ and error in days.}
  \label{sfig:corr-day}
\end{minipage}
\hfill
\begin{minipage}[t]{0.45\textwidth}
  \centering
  \includegraphics[width=\textwidth]{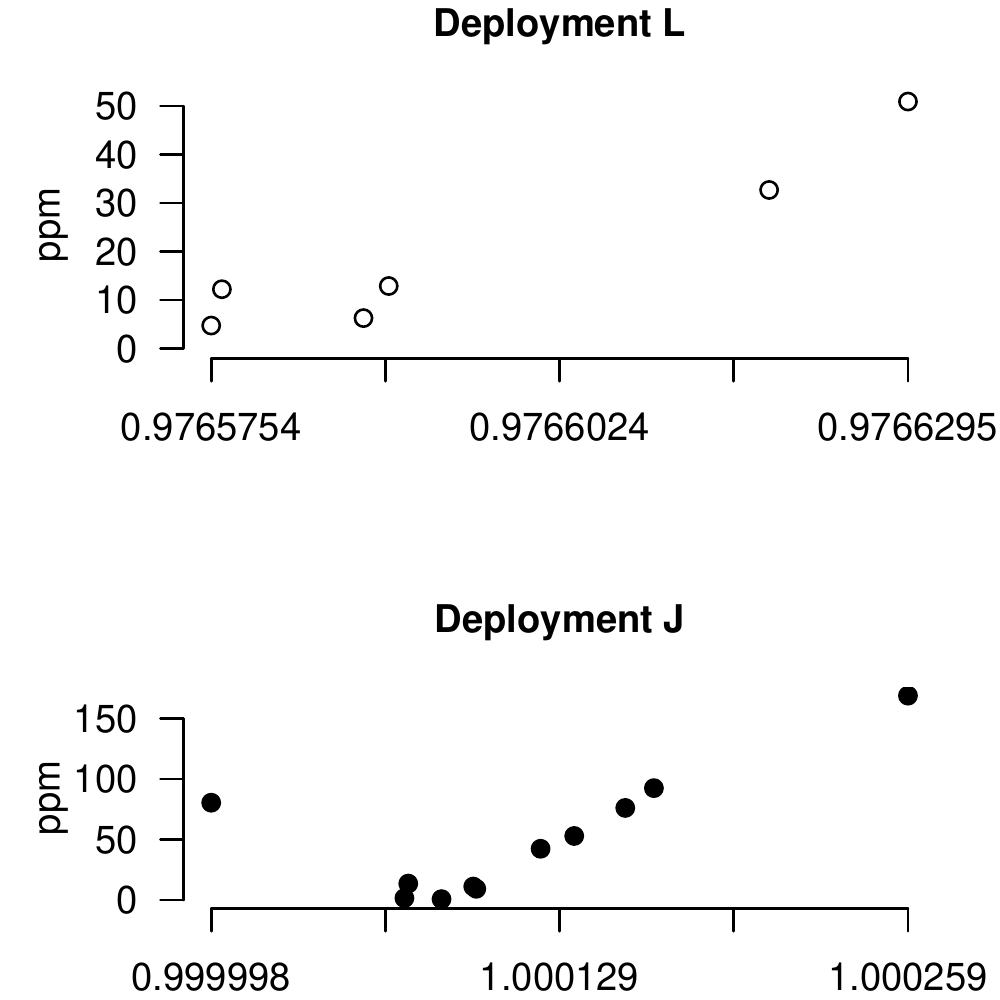}
  \caption{$\alpha$ estimates from Sundial and estimation errors in ppm.}
  \label{sfig:alpha}
\end{minipage}
\end{figure}

\subsection{Reconstructing Global Timestamps using Sundial}

We evaluate Sundial using data from the segments shown in
Figure~\ref{fig:deploy}. Specifically, we evaluate the accuracy of
the timestamps reconstructed by Sundial as though the start time of
these segment is unknown (similar to the case of a mote reboot) and no
anchor points are available. 
Since we make no assumptions of the segment start-time, a very
large model (solar) signal needs to be considered to find the correct
shift (lag) for the day alignment.

\subsubsection{Evaluation Metrics:} We divide the timestamp reconstruction
error to: (a) error in days; and (b) error in minutes within the
day. The error in minutes is computed as the root mean square error
($RMSE_{min}$) over all the measurements. We divide the reconstruction
error into these two components, because this decoupling
naturally reflects the accuracy of estimating the $\alpha$ and $\beta$
parameters. Specifically, if the $\alpha$ estimate were inaccurate,
then, as Figure~\ref{fig:drift} suggests, the reconstruction error
would grow as a function of time. In turn, this would result in a
large root mean squared error in minutes within the day over all the
measurements. On the other hand, a low $RMSE_{min}$ corresponds to an
accurate estimate for $\alpha$. Likewise, inaccuracies in the
estimation of $\beta$ would result in large error in days.

\subsubsection{Results:} Figures~\ref{sfig:dayerr} and~\ref{sfig:minerr}
summarize Sundial's accuracy results. Overall, we find that longer
segments show a lower day error. Segments belonging to the {\em L}
deployment span well over a year and the minimum day error is 0 while
the maximum day error is 6. In contrast, most of the segments for
deployment {\em J} are less than 6 months long and the error in days
for all but two of those segments is less than one week.
Figure~\ref{sfig:corr-day} presents the relationship between the
maximum correlation ($\rho_{max}$) and the day error.  As $\rho_{max}$
measures how well we are able to match the LOD pattern for a node with
the solar LOD pattern, it is not surprising that high correlation is
generally associated with low reconstruction error.  The $RMSE_{min}$
obtained for each of the segments in deployment {\em L} is very low
(see Figure~\ref{sfig:minerr}) . Remarkably, we are able to achieve an
accuracy ($RMSE_{min}$) of under a minute for the majority of the
nodes of the {\em L} deployment even though we are limited by our
sampling frequency of 20 minutes. Moreover, $RMSE_{min}$ error is
always within one sample period for all but one segment.

Interestingly, we found that the $\alpha$ values for the two
deployments were significantly different. This disparity can be
attributed to differences in node types and thus clock
logic. Nonetheless, Sundial accurately determined $\alpha$ in both
cases. Figure~\ref{sfig:alpha} presents the $\alpha$ values for the
two deployments. We also show the error between the $\alpha$ obtained
using Sundial and the $\alpha$ value obtained by fitting the good
anchor points sampled by the gateway (i.e., ground truth fit).  The
ppm error for both the deployments is remarkably low and close to the
operating error of the quartz crystal.

\subsection{Impact of Segment Length}
\label{subsec:seg-len}

\begin{figure}[t]
\begin{minipage}[t]{0.45\textwidth}
  \centering
  \includegraphics[width=\textwidth]{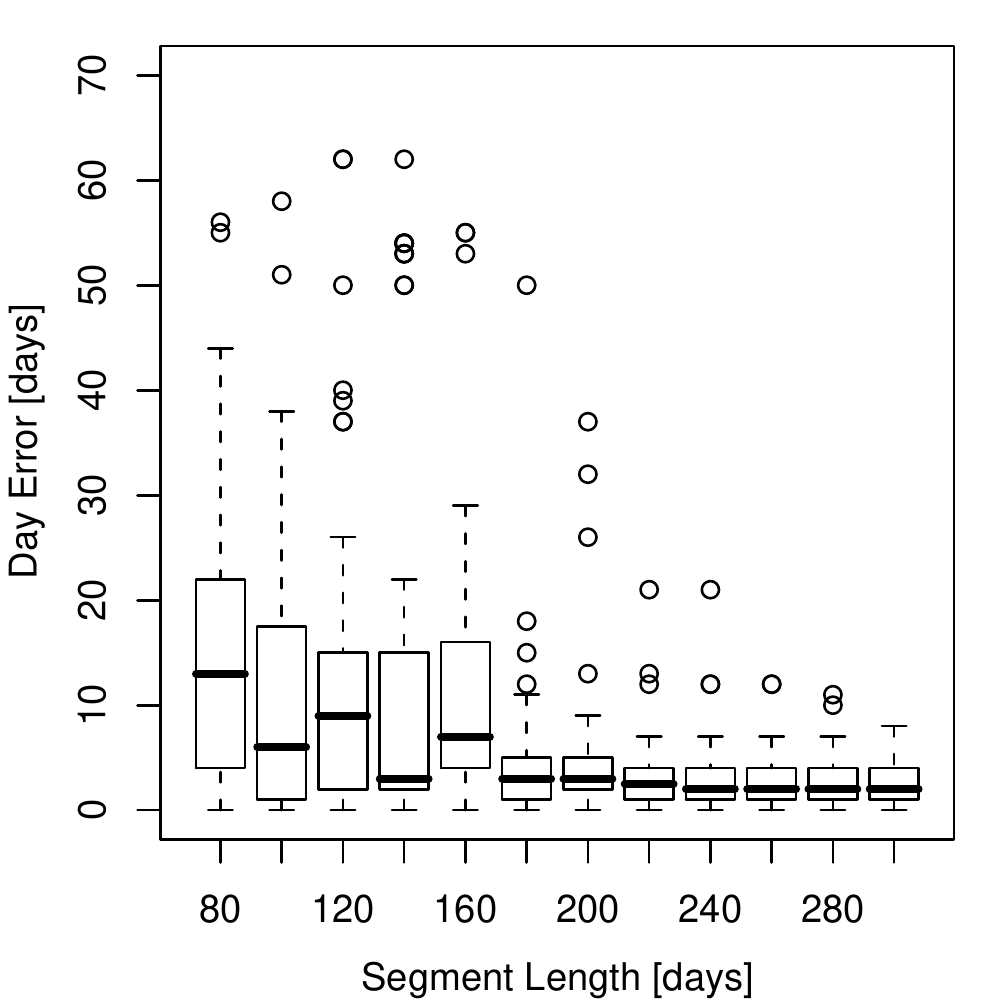}
  \caption{Error in days as a function of segment size.}
  \label{sfig:vdays-dayerr}
\end{minipage}
\hfill
\begin{minipage}[t]{0.45\textwidth}
  \centering
  \includegraphics[width=\textwidth]{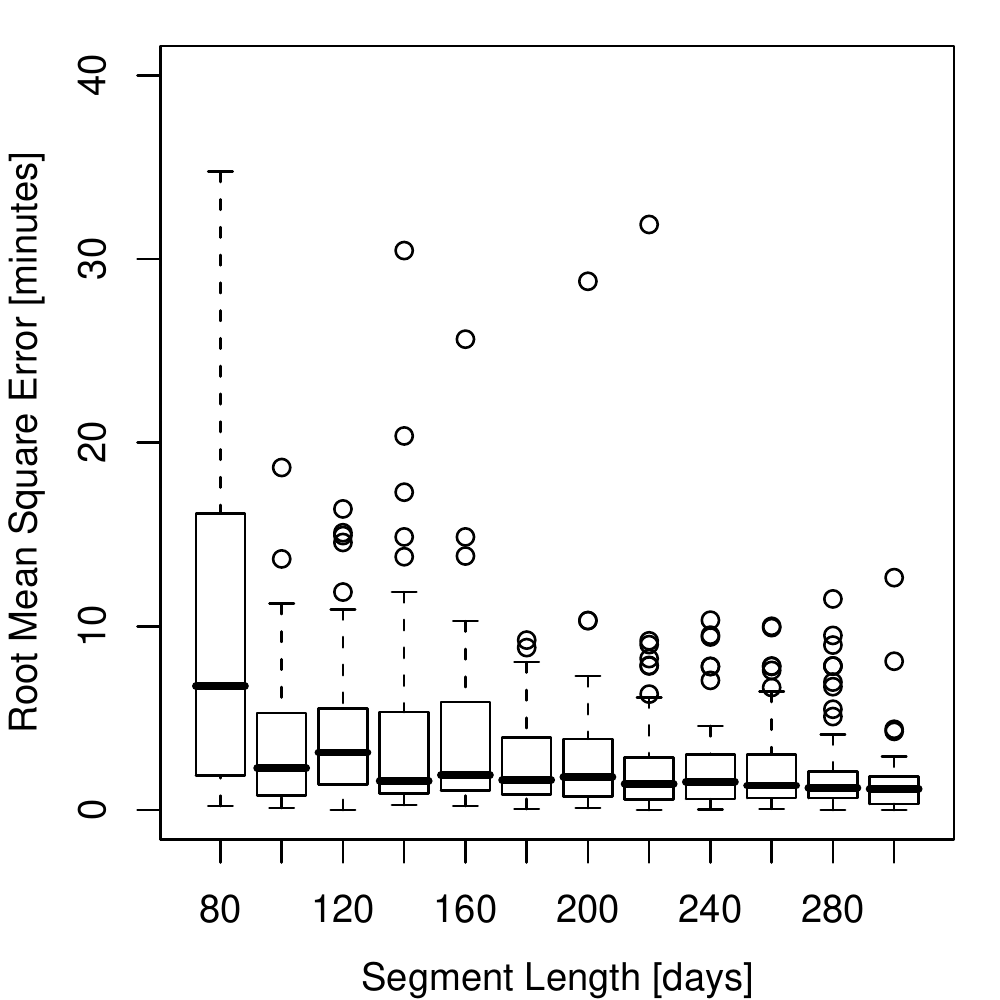}
  \caption{Error in minutes ($RMSE_{min}$) as a function of segment size.}
  \label{sfig:vdays-minerr}
\end{minipage}
\end{figure}

Sundial relies on matching solar patterns to the ones observed by the
light sensors. The natural question to ask is: what effect does the
length of segment have on the reconstruction error. We address this
question by experimenting with the length of segments and observing
the reconstruction error in days and $RMSE_{min}$.  We selected data
from three long segments from deployment {\em L}. To eliminate bias,
the start of each shortened segment was chosen from a uniform random
distribution. Figure~\ref{sfig:vdays-minerr} shows that the
$RMSE_{min}$ tends to be remarkably stable even for short
segments. One concludes that even for short segment lengths, Sundial
estimates the clock drift ($\alpha$)
accurately. Figure~\ref{sfig:vdays-dayerr} shows the effect of segment
size on day error. In general, the day error decreases as the segment
size increases. Moreover, for segments less than 150 days long, the
error tends to vary considerably.

\subsection{Day Correction}

\begin{figure}[t]
  \centering
  \includegraphics[scale=0.5]{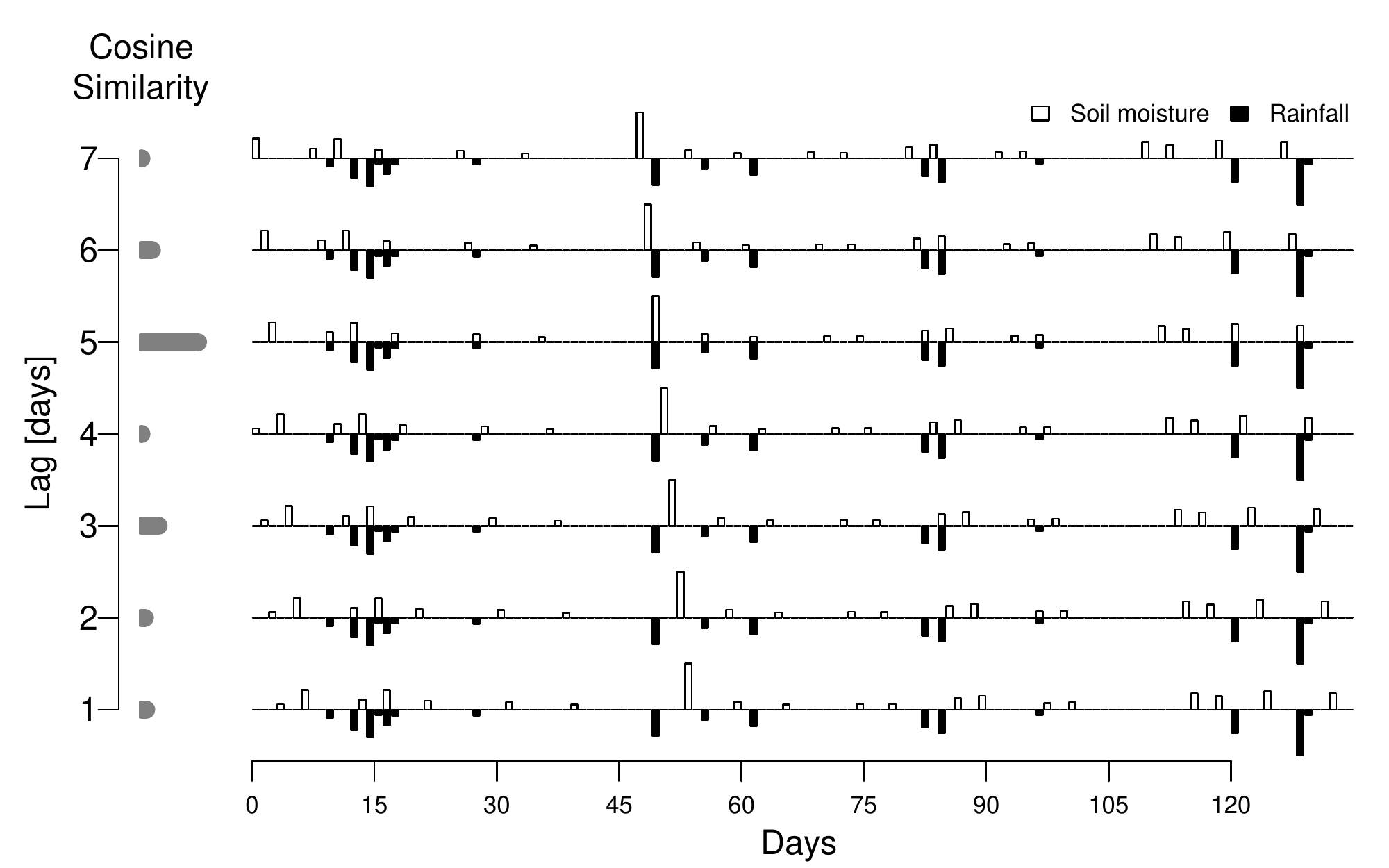}
  \caption{An illustration of the cosine similarity
    ($\theta_{SM-PPT}$) values for seven different day lags between
    moisture and rainfall vectors. $\theta_{SM-PPT}$ peaks at the
    correct lag of five days, providing the correct day adjustment.}    
  \label{fig:rainfix}
\end{figure}

The results so far show that 88\% (15 out of 17) of the motes have a
day offset of less than a week. Next, we demonstrate how global events
can be used to correct for the day offset. We looked at soil moisture
data from eight motes of the {\em J} deployment after obtaining the
best possible timestamp reconstruction. Specifically, we correlated
the motes' soil moisture data with rainfall data to correct for the
day offset. We used rainfall data from a period of 133 days, starting
from December 4, 2007, during which 21 major rain events occurred.  To
calculate the correlation, we created weighted daily vectors for soil
moisture measurements ($SM$) whose value was greater than a certain
threshold and similarly rainfall vectors having a daily precipitation
($PPT$) value of greater than 4.0 cm. Next, we extracted the lag at
which the cosine angle between the two vectors (cosine similarity,
$\theta_{SM-PPT}$) is maximum. This method is inspired by the
well-known document clustering model used in the information retrieval
community~\cite{SALTON}. Note that we computed $\theta_{SM-PPT}$ for a
two-week window ($\pm$ seven days) of lags and found that seven out of
the eight motes could be aligned perfectly.  Figure~\ref{fig:rainfix}
illustrates the soil moisture vectors, rainfall vectors and the
associated $\theta_{SM-PPT}$ for seven lags for one of the
segments. Note that $\theta_{SM-PPT}$ peaks at the correct lag of
five, leading to the precise day correction.  While we use soil
moisture to illustrate how global events can be used to achieve
macro-level clock adjustments, other modalities can also be used based
on the application's parameters.

%% file: relwork.tex
\section{Related Work}
\label{sec:relwork}

This study proposes a solution to the problem of postmortem timestamp
reconstruction for sensor measurements. To our knowledge, there is
little previous work that addresses this problem for deployments that
span a year or longer. Deployment length can be an issue because the
reconstruction error monotonically increases as a function of time
(cf. Sec.\ref{subsec:reboot}). The timestamp reconstruction problem
was first introduced by Werner-Allen et al. who provided a detailed
account of the challenges they faced in synchronizing mote clocks
during a 19-day deployment at an active
volcano~\cite{GWA06}. Specifically, while the system employed the FTSP
protocol to synchronize the network's motes, unexpected faults forced
the authors to rely on an offline {\em time rectification} algorithm
to reconstruct global timestamps.

While experiences such as the one reported in~\cite{GWA06} provide
motivation for an independent time reconstruction mechanism such as
the one proposed in this paper, the problem addressed by Werner-Allen
et al. is different from the one we aim to solve. Specifically, the
volcano deployment had access to precise global timestamps (through a
GPS receiver deployed at the site) and used linear regression to
translate local timestamps to global time, once timestamp outliers
were removed. While RGTR can also be used for outlier detection and
timestamp reconstruction, Sundial aims to recover timestamps in
situations where a reliable global clock source is not available.

Finally, Chang et. al.~\cite{CHANG08} describe their experiences with
motes rebooting and resetting of logical clocks, but do not furnish
any details of how they reconstructed the global timestamps when this
happens.

%% file: concl.tex
\section{Conclusion}
\label{sec:concl}

In this paper we present Sundial, a method that uses light sensors to
reconstruct global timestamps. Specifically, Sundial uses light
intensity measurements, collected by the motes' on-board sensors, to
reconstruct the length of day (LOD) and noon time throughout the
deployment period. It then calculates the slope and the offset by
maximizing the correlation between the measurement-derived LOD series
and the one provided by astronomy. Sundial operates in the absence of
global clocks and allows for random node reboots. These features make
Sundial very attractive for environmental monitoring networks deployed
in harsh environments, where they operate disconnected over long
periods of time. Furthermore, Sundial can be used as an independent
verification technique along with any other time reconstruction
algorithm.

Using data collected by two network deployments spanning a total of
2.5 years we show that Sundial can achieve accuracy in the order of a
few minutes. Furthermore, we show that one can use other global events
such as rain events to correct any day offsets that might exist. 
As expected, Sundial's accuracy is closely related to the segment
size. In this study, we perform only a preliminary investigation on
how the length of the segment affects accuracy. An interesting
research direction we would like to pursue is to study the
applicability of Sundial to different deployments. Specifically, we
are interested in understanding how sampling frequency, segment
length, latitude and season (time of year) collectively affect
reconstruction accuracy.

Sundial exploits the correlation between the well-understood solar
model and the measurements obtained from inexpensive light sensors. In
principle, any modality having a well-understood model can be used as
a replacement for Sundial. In the absence of a model, one can exploit
correlation from a trusted data source to achieve reconstruction,
e.g., correlating the ambient temperature measurement between the motes
with data obtained from a nearby weather station.  However, we note
that many modalities (such as ambient temperature) can be highly
susceptible to micro-climate effects and exhibit a high degree a
spatial and temporal variation. Thus, the micro-climate invariant
solar model makes light a robust modality to reconstruct timestamps in
the absence of any sampled anchor points. 

Finally, we would like to emphasize the observation that most
environmental modalities are affected by the diurnal and annual solar
cycles and not by the human-created universal time. In this regard,
the time base that Sundial establishes offers a more natural reference
basis for environmental measurements.

%% file: ack.tex
\section*{Acknowledgments}
\label{sec:ack}

We would like to thank Yulia Savva (JHU, Department of Earth and
Planetary Science) for helping us identify the timestamp
reconstruction problem.  This research was supported in part by NSF
grants CNS-0546648, CSR-0720730, and DBI-0754782. Any opinions,
finding, conclusions or recommendations expressed in this publication
are those of the authors and do not represent the policy or position
of the NSF.